\documentclass[12pt]{amsart}
\usepackage{fullpage,amssymb,epic,eepic}
\def\draft{n}

\theoremstyle{plain}

\newtheorem{theorem}{Theorem}
\newtheorem{proposition}{Proposition}[section]
\newtheorem{lemma}[proposition]{Lemma}
\newtheorem{corollary}[proposition]{Corollary}

\theoremstyle{definition}
\newtheorem{definition}[proposition]{Definition}

\theoremstyle{remark}

\newtheorem{remark}[proposition]{Remark}

\def\printname#1{
	\if\draft y
		\smash{\makebox[0pt]{\hspace{-0.5in}
			\raisebox{8pt}{\tt\tiny #1}}}
	\fi
}

\newlength{\standardunitlength}
\catcode`\@=11
\long\def\@makecaption#1#2{%
    \vskip 10pt
    \setbox\@tempboxa\hbox{
      \small\sf{\bfcaptionfont #1. }\ignorespaces #2}%
    \ifdim \wd\@tempboxa >\captionwidth {%
        \rightskip=\@captionmargin\leftskip=\@captionmargin
        \unhbox\@tempboxa\par}%
      \else
        \hbox to\hsize{\hfil\box\@tempboxa\hfil}%
    \fi}
\font\bfcaptionfont=cmssbx10 scaled \magstephalf
\newdimen\@captionmargin\@captionmargin=2\parindent
\newdimen\captionwidth\captionwidth=\hsize
\catcode`\@=12

\def\lbl#1{\label{#1}\printname{#1}}

\def\eqdef{\overset{\text{def}}{=}}

\def\eqdef{\overset{\text{def}}{=}}

\def\ds{\displaystyle}

\def\BZ{\mathbb Z}
\def\BQ{\mathbb Q}
\def\BR{\mathbb R}
\def\BC{\mathbb C}

\def\C{\mathcal C}

\def\s{\sigma}
\def\pa{\partial}
\def\p{\partial}
\def\ep{\epsilon}
\def\z{\zeta}
\def\l{\lambda}
\def\w{\omega}
\def\th{\Theta}

\def\hB{\widehat{\text{Ber}}}

\def\ft{\mathfrak t}

\def\fs{\mathfrak s}
\def\di{\diamond}

\def\b{\beta}
\def\cs{\check{\sigma}}
\def\Ber{\text{Ber}}

\def\pd#1{\frac{\partial }{\partial #1}}

\def\tByt#1#2#3#4{ \begin{pmatrix}
                   #1 & #2 \\
                   #3 & #4 
                   \end{pmatrix} }

\def\2by1#1#2{ \bigl(\begin{smallmatrix}
                   #1  \\
                   #2 
                   \end{smallmatrix} \bigr)}

\def\tByo#1#2{ \begin{pmatrix}
                   #1  \\
                   #2 
                   \end{pmatrix} }

\renewcommand{\)}{)\!)}
\renewcommand{\]}{{]\!]}}
\renewcommand{\[}{{[\![}}

\def\lp{\langle}
\def\rp{\rangle}

\def\num{\gamma}
\def\denom{\delta}
\def\zz{\mathbb Z}
\def\conevar{\varpi}

\begin{document}


\title[Values of zeta functions, Dedekind sums and toric geometry]{Values of 
      zeta 
      functions at negative  integers, Dedekind sums and  toric 
       geometry}

\author{Stavros Garoufalidis}
\address{Department of Mathematics \\
         Harvard University \\
         Cambridge, MA 02138, USA. }
\email{stavros@math.harvard.edu}
\thanks{The  authors were partially supported by NSF grants 
       DMS-95-05105 and DMS-95-08972, respectively.\newline
       This and related preprints can also be obtained at
{\tt http://www.math.brown.edu/$\sim$stavrosg } \newline
1991 {\em Mathematics Classification.} Primary 11M06. Secondary 14M25,
11F20.
}
\author{James E. Pommersheim}
\address{Department of Mathematical Sciences\\
        New Mexico State University\\
        Las Cruces, NM 88003, USA. }
\email{jamie@math.nmsu.edu}

\dedicatory{Dedicated to our teacher, W. Fulton.}

\date{This edition: September 12, 1998; First edition: November 22, 1997 }

\begin{abstract}
We study relations among   special values of  zeta functions, invariants
of toric varieties, and  generalized Dedekind sums.  In particular,
we use invariants arising in the Todd class of a toric variety to
give a new explicit formula for the values of the zeta function of
a real quadratic field at nonpositive integers.  We also
express these invariants in terms of the generalized Dedekind
sums studied previously by several authors.
The paper includes conceptual proofs of these relations
and explicit computations of the various zeta values and Dedekind sums 
involved.
\end{abstract}

\maketitle

\tableofcontents


\section{Introduction}
\lbl{sec.intro} 

In the  present paper, we study relations among
 special values of  zeta functions of real quadratic fields, 
properties of generalized Dedekind sums and   Todd classes of toric 
varieties. The main theme of the paper is
the use of toric geometry to explain in a conceptual way
properties of the values of
zeta functions and Dedekind sums, as well to provide explicit computations.

Both toric varieties and zeta functions associate numerical invariants
to cones in lattices;
however, with different motivations and applications.
Though we will focus on the case of two-dimensional
cones in the present paper, we introduce notation and definitions that are 
valid for
cones of arbitrary dimension. The reasons for this added generality is clarity,
as well as preparation for the results of a subsequent publication.

\subsection{Zeta functions}
\lbl{sub.history}



We begin by  reviewing the first source of numerical invariants of cones:
the study of zeta functions.

Given a number field $K$,  its  zeta function is defined (for $Re(s)$ 
sufficiently large) 
by:
\begin{equation*}
\z(K,s)=\sum_{\alpha} \frac{1}{Q(\alpha)^s}
\end{equation*}
where the summation is over all nonzero ideals,
 and $Q$ is the norm. The above function
admits a meromorphic continuation in $\BC$, with a simple pole
at $s=1$, and regular everywhere else.
Lichtenbaum \cite{Li}
conjectured a specific behavior of the zeta function 
$\z(K,s)$ at nonpositive
integers related to the global arithmetic of
the number field. In the special 
 case  of a totally real number field $K$,
Lichtenbaum conjectured that the values of the zeta function at negative
integers are rational numbers which involve the rank of the algebraic 
(or \'etale) $K$-theory of $K$. It turns out  (see 
\cite{Sh1} and \cite[Section 2]{Za1}) that 
for a totally real field $K$ the zeta function
can be decomposed as a  sum $\z(K,s)= \sum_I \z(I,s)$ where the sum 
is over the (finite set of) narrow ideal classes. For each narrow ideal
class $I$, there is a finite set  $\{ \s^{\circ}_{i,I} \}_i$ of open cones
in a rank $[K:\BQ]$
lattice $M_I$ in $K$ such that $\z(I^{-1},s)=\sum_i \z_{Q,\s^{\circ}_{i,I}}(s)$
where $Q$ is a multiple of the norm  on $K$ and where for an open
cone $\tau^{\circ}$ we have
$$
\z_{Q,\tau^{\circ}}(s)=\sum_{a \in \tau^\circ \cap M_I}\frac{1}{Q(a)^s}.
$$
The problem of calculating the zeta values $\z_{Q,\tau^{\circ}}(-n)$ for
$n \geq 0$ for all triples $(M,Q,\tau^{\circ})$ that arise from totally real
fields has attracted a lot of attention by several authors.
Klingen \cite{Kl} and Siegel \cite{Si1,Si2}
 (using analytic methods) proved that the values of the zeta functions
of totally real fields 
at nonpositive integers are rational numbers and provided  an algorithm
 for calculating  them. Meanwhile, 
Shintani \cite{Sh1} (using algebraic and combinatorial methods)
gave an independent calculation of the zeta values $\z_{Q,\tau^{\circ}}(-n)$
for an arbitrary $Q$ which is a product of linear forms.
Meyer \cite{Me} and Zagier \cite{Za3} gave another calculation
of  the zeta values $\z_{Q,\tau^{\circ}}(-n)$ for rank two lattices. 
Related results have also been obtained by P. Cassou-Nogu\`{e}s, \cite{C-N1,
C-N2, C-N3}. More recently,
Hayes \cite{Hs}, Sczech \cite{Sc1, Sc2} and Stevens \cite{St}
have constructed $\text{PGL}_m(\BQ)$ cocycles which, among other
things, provide a calculation of the zeta values (at nonpositive integers)
of totally real number fields in terms of generalized Dedekind sums.

We now specialize to  the case of
real quadratic fields.  Given a narrow ideal class $I$, there is
a  cone  $\s_I$ (i.e., the convex hull of two rays from the origin) 
in a rank two lattice $M \subseteq K$, such that:
$\z(I^{-1},s)= \z_{Q,\s_I}(s).$ 
Here $Q$ is a multiple of the norm form and $\z_{Q,\s}$ is given
by the following definition:

\begin{definition}
\lbl{def.zeta}
Let $M$ be an two-dimensional lattice, and let $V$ denote the real
vector space $ M_{\BR}=M\otimes \BR$. Let
 $Q: M_{\BR} \to \BR$
 be  a {\em nonzero quadratic homogenous function} 
i.e., a function satisfying
$Q(a v)=a^2 Q(v)$ for $a \in \BR , 
v \in M_{\BR}$ and such that $Q$ is not identically $0$.
We set
\begin{equation}
\z_{Q,\tau}(s)= \sum_{ a \in \tau \cap M} \frac{wt(\tau, a)}{Q(a)^s}
\end{equation}
where $\tau$ is a rational two-dimensional cone in $M$ 
and $wt(\tau, \cdot ): M
\to \BQ$ is the weight function defined by:
\begin{equation}
\lbl{eq.we}
wt(\tau, a)=
\begin{cases}
1 & \text{if $a$ lies in the interior of $\tau$} \\
1/2 & \text{if $a$ lies in the boundary of $\tau$, and  $a \neq 0$} \\
0 & \text{otherwise.}
\end{cases}
\end{equation}
\end{definition}
\noindent
We will assume throughout that $Q$ is positive on all of $\tau$;
that is, for all
$a\in\tau$, $a\neq 0$, we have $Q(a)>0$.
The above zeta function, 
defined for $Re(s)$ sufficiently large, \cite{Za2} can be 
analytically continued to a meromorphic function on $\BC$, regular everywhere
except at $1$. 

In \cite{Za2}, Zagier showed that
after possibly multiplying $M$ by a totally positive number and
$Q$ by a nonzero rational number (which only multiplies values of the zeta
function by a nonzero constant),  every triple
$(M,Q,\tau)$ which arise from real quadratic  
field can be constructed  by means of
 a finite sequence $b=(b_0, \dots,b_{r-1})$
of integers greater than $1$ and not all equal to $2$,
which will be denoted by $(M_b,Q_b,\tau_b)$. In addition,
there are canonical vectors $A_0, \dots, A_r$ in $M_b$ (that depend on $b$)
that 
subdivide the cone $\tau_b= \lp A_0, A_r \rp$ (i.e., the cone whose extreme
rays are $A_0$ and $A_r$) in $M_b$ in $r$ nonsingular cones
$\lp A_i, A_{i+1} \rp$. For a detailed discussion, see
Section \ref{sub.real}.

We will now give an explicit formula for the zeta values 
$ \z_{Q_{b}, \tau_b}(-n)$ (for $ n \geq 0$), which in turn enables one
to calculate the values  of the zeta function
of a real quadratic field at nonpositive integers.  While both the motivation
and the proof involve concepts from the theory of toric geometry, the formula
can be stated and understood without a knowledge of toric varieties.  We will
state the formula here, and in the next subsection, we will discuss concepts
from toric geometry which are necessary for the proof and which lead
to a conceptual understanding of the present formula.

Let $\l_m$ be defined by the power series:
\begin{equation}
\lbl{eq.lm}
\frac{h}{1-e^{-h}} = \sum_{m=0}^\infty \l_m h^m
\end{equation}
thus we have: $\l_m = (-1)^m B_m/m!$ where $B_m$ is the $m^{th}$
 Bernoulli number. (See also Definition \ref{def.dedsum} below.)  Note
that if $m>1$ is odd, then $\lambda_m=0$.
For $d\geq 2$ even, define homogeneous polynomials $P_d(X,Y), R_d(X,Y)$ 
of degree
$d-2$ by:
\begin{eqnarray*}
P_d(X,Y)& = & \epsilon_{d}\sum_{i+j=d,i,j>0} \lambda_i \lambda_j X^{i-1} Y^{j-1},\\
R_d(X,Y)& = & \frac{X^{d-1}+Y^{d-1}}{X+Y}=X^{d-2}-X^{d-3}Y+\dots+Y^{d-2},
\end{eqnarray*}
\noindent
where $\epsilon_{2}=-1$ and $\epsilon_{d}=1$ for all even $d>2$.

We then have:
\begin{theorem}
\lbl{cor.zv}
The values  $\z_{Q_{b}, \tau_b}(-n)$ for $n \geq 0$ are polynomials 
in  $b_i$ with rational coefficients 
(symmetric under cyclic permutation of the $b_i$) 
 given explicitly as follows:
\begin{eqnarray}
\lbl{eq.zva}
 \z_{Q_{b}, \tau_b}(-n) & = &
 (-1)^n n! \left\{
\sum_{i=1}^r P_{2n+2}\left(\pd x , \pd y \right) 
\di e^{- Q_b( x A_{i-1} + y A_i)}
\right. \\
& & \left.
+ \l_{2n+2} \sum_{i=1}^r b_i R_{2n+2}\left(\pd x , \pd y \right) \di   
e^{- Q_b( x A_{i-1} + y A_{i+1})} \right\} 
\end{eqnarray}
The diamond symbol indicates the evaluation of a partial derivative at the 
origin.
\end{theorem}
\noindent
In particular, we obtain the formula of Zagier \cite[Equation 3.3]{Za4}:
\begin{eqnarray}
\lbl{eq.z0}
\z_{Q_{b}, \tau_b}(0) & = & \frac{1}{12} \sum_{i=0}^{r-1} (b_i-3).
\end{eqnarray}

\subsection{Toric geometry}
\lbl{sub.toric}

The second source of numerical invariants of cones
comes from the theory of {\em toric varieties}. 
For a general reference on toric varieties, see \cite{Dan} or \cite{Fu}.

Founded in
the 1970s, the subject of toric
varieties provides a strong link between algebraic geometry and the theory of 
convex bodies in a lattice.  To each lattice polytope (the convex hull of
a finite set of lattice points) is associated an algebraic variety
with a natural torus action.  This correspondence enables one to translate
important properties and theorems about lattice polytopes into the language
of algebraic geometry, and vice-versa. One important example of this
is the very classical problem of counting the number of lattice points in a 
polytope.  The early pioneers in the subject of toric varieties found that
this problem, translated into algebraic geometry, becomes the problem 
of finding the Todd class of a toric variety.
Much progress  has been made over the past ten years on the Todd
class problem,  and this has led to a greater understanding of the 
lattice point counting question.

One approach to computing the Todd class of a toric variety
is the fundamental work of R.
Morelli \cite{Mo}.  He settled a question of Danilov by proving a local
formula expressing the Todd class of a toric variety as a cycle.  Another
approach, introduced in by the first author in \cite{P1, P2}, is to express 
the Todd class as a
polynomial in the torus-invariant cycles.  Dedekind sums appear 
as coefficients in these polynomials, and this leads to lattice point
formulas in terms of Dedekind sums, as well as new reciprocity laws
for Dedekind sums \cite{P1}. Cappell and Shaneson \cite{CS} subsequently
announced an extension of the program of \cite{P1} in which they proposed 
formulas for the Todd class of a toric variety in all dimensions.  
In \cite{P3} it was shown that the
polynomials of \cite{P1, P2} 
can be expressed nicely as the truncation of a certain power
series whose coefficients were shown to be polynomial-time
computable using an idea of Barvinok \cite{Ba}. A beautiful power series 
expression for the equivariant Todd
class of a toric variety was given by Brion and Vergne
in \cite{BV2}.  Guillemin \cite{Gu} also proved similar Todd class formulas
from a symplectic geometry point of view.  Furthermore, in \cite{BV}, Brion
and Vergne use the Todd power series of \cite{BV2} to give a formula
for summing any polynomial function over a polytope. 

 The power
series studied in \cite{P3} and \cite{BV,BV2} are, in fact, identical
(see Section \ref{sec.thm.1}) and play a central role in the present
paper.  A detailed discussion of the properties of these power
series, which we call the {\em Todd power series of
a cone}, is contained in Section \ref{sec.thm.1}.  In this section, we state
two theorems about the Todd power series of a two-dimensional cone
that we will need in our study of the zeta function.

Given independent rays $\rho_1, \dots , \rho_n$ from the origin
in an $n$-dimensional lattice
$N$, the convex hull of the rays in the vector space $N_{\BR}
=N\otimes\BR$
forms an $n$-dimensional cone $\s= \lp \rho_1, \dots , \rho_n \rp $. Cones of
this type (that is, ones generated by linearly independent rays)
are called {\em simplicial}. Let $\C^n(N)$ denote the set of $n$-dimensional
simplicial cones of $N$ with ordered rays. There is then a canonical
function 
\begin{equation*}
\ft : \C^n (N) \to \BQ\[x_1, \dots, x_n \],
\end{equation*}
invariant under lattice automorphisms, which associates to each cone $\sigma$ a power series $\ft_{\sigma}$ with rational
coefficients, called the Todd power series of $\sigma$.
Several ways of characterizing this function are given in Section
\ref{sec.thm.1}.  These include an $N$-additivity property (See Proposition
\ref{prop.j1}), an exponential sum over the cone (Proposition \ref{prop.j2})
and an explicit cyclotomic sum formula (Proposition \ref{prop.j3}.)

In the case of a two-dimensional cone $\s$, the coefficient of $xy$ of
$\ft_\s(x,y)$ was identified as a {\em Dedekind sum} \cite{P1}.
Furthermore, in \cite{P1} it was shown that
{\em reciprocity formulas} for Dedekind sums follow from 
an $N$-additivity formula for $\ft$.  Zagier's higher-dimensional 
Dedekind sums \cite{Za4} were later shown to appear as coefficients \cite{BV2}.
  It is natural to conjecture that, in the two-dimensional case,
{\em all} the coefficients
of $\ft$ are generalized Dedekind sums and
that reciprocity properties of generalized Dedekind sums will be related
to the $N$-additivity formula of $\ft$.
Indeed, this is the case:  see Theorem \ref{thm.f=s}.

We now present an explicit link between Todd power series and
zeta functions.  The following theorem expresses the values of zeta
function at negative integers in terms of the Todd power series of a 
two-dimensional cone.  Note that the idea of considering the Todd power 
series as
a differential operator, and applying this to an integral over a shifted 
cone is not new: this was introduced in \cite{KP} and developed further
in \cite{BV}.

First we introduce some notation which is standard in the theory of
toric varieties.  If $\tau$ is a cone in a lattice
$M$, the dual cone $\check\tau$ is a cone in the dual lattice 
$N=\text{Hom}(M,\BZ)$,
defined by
$$
\check\tau=\{v\in N | \lp v, u \rp \ge 0\ \text{for all}\ u\in\tau\}.
$$
The dual of an $n$-dimensional simplicial cone $\tau=\lp\rho_1,\dots,
\rho_n\rp$ is generated by the rays $u_i$
a normal to the $n-1$-dimensional faces of $\tau$; hence $\check\tau$
is also a simplicial cone. Given $h=(h_1,\dots,h_n)\in\BR^n$, we define
the {\em shifted cone} $\tau(h)$ to be the following cone in $M_{\BR}$:
$$\tau(h)=
\{ m \in M_{\BR}| \lp u_i, m \rp \geq -h_i \text{ for all }
i=1, \dots, n \}.
$$
Here (and throughout) we have identified each ray
$u_i$ with the primitive lattice point on that ray, that is, the
nonzero lattice point on $u_i$ closest to the origin.

Given an $n$-dimensional simplicial cone $\tau=\lp \rho_1,\dots,\rho_n\rp$ 
in an 
$n$-dimensional lattice $M$, the
{\em multiplicity} of $\tau$, denoted by $\text{mult}(\tau)$,
is defined to be the index in $M$ of the sublattice $\BZ\rho_1+\dots
+\BZ\rho_n$ (again identifying the rays $\rho_i$ with their 
primitive lattice points.) Thus
the multiplicity  of $\tau$ is simply the volume of the parallelepiped formed
by the vectors from the origin to the primitive lattice points
on the rays of $\tau$.

\begin{theorem}
\lbl{thm.1}
Let $\tau$ be a two-dimensional cone of multiplicity $q$
 in a two-dimensional lattice
$M$, and let $\sigma=\check\tau$ be the dual cone in 
$N=Hom(M,\BZ)$.
Then for all $n \geq 0$,
we have:
$$\z_{Q,\tau}(-n)=
(-1)^n n!
\left\{ (\ft_{\s})_{(2n+2)} \left(\frac{\p}{\p h_1},\frac{\p}{\p h_2} \right) 
-\delta_{n,0}
\frac{q}{2} \frac{\p^2}{\p h_1 \p h_2} \right\} 
\di  \int_{\tau(h)} exp(-Q(u))du,
$$
where all derivatives above are evaluated at $h_1=h_2=0$
and $\delta_{n,0}=1$ (resp. $0$) if $n=0$ (resp. $ n \neq 0$).
\end{theorem}
\noindent
In this equation, $(\ft_{\s})_d$ denotes the degree $d$ part of  
 the Todd power series thought of as an (infinite order)
 constant coefficients differential operator acting on the function
$h \to \int_{\tau(h)} exp(-Q (u))du$.

The coefficients of $\ft_{\s}$ for a two-dimensional cone may be expressed
explicitly in terms of continued fractions.  We now state this formula.

Let $\sigma$ be a two-dimensional cone in a lattice $N$.
Then there exist a (unique) pair of relatively
prime integers $p,q$ with $0< p \le q$ such
that $\sigma$ is lattice equivalent (equivalent under the automorphism group
of the lattice) to the cone $\sigma_{(p,q)}
\eqdef \lp (1,0),(p,q) \rp$ in 
$\BZ^2$. Such a cone will be called a cone of type $(p,q)$. Let
$ b_i,  h_i,  k_i,  X_i$  are defined in terms of  the negative-regular
continued
fraction expansions:
\begin{equation}
\lbl{eq.bis}
\frac{q}{p} = [b_1,\dots,b_{r-1}], \qquad
\frac{ h_i}{ k_i}  \eqdef  [b_1,\dots,b_{i-1}], \qquad
 X_i  \eqdef  -  h_i x+(q  k_i -p  h_i)y.
\end{equation}
(Throughout, we use
such bracketed lists to denote finite negative-regular continued fractions.)
We then have the following continued fraction expression for the degree $d$  
part $(\ft_{\sigma})_d$  of the Todd power series $\ft_{\sigma}$.

\begin{theorem}
\lbl{thm.toddd}
For $\s$ of type $(p,q)$ as above, and for  $d\geq 
2$  an even integer, we have:
\begin{eqnarray*}
(\ft_{\sigma})_{d}(x,y)& = & qxy \sum_{i=1}^r
P_d(  X_{i-1}, X_i) +\lambda_d qxy\sum_{i=1}^{r-1} b_i R_d( X_{i-1},
 X_{i+1}) \\ & &
-\lambda_d(x  X_1^{d-1}+y  X_{r-1}^{d-1})+\frac12 \delta_{d,2}qxy.
\end{eqnarray*}
If $d\geq 1$ is odd, then $ (\ft_{\s})_{d}(x,y)
=\frac12\lambda_{d-1}q^{d-2}xy(x^{d-2}+y^{d-2}). $
\end{theorem}

\begin{remark}
We have stated the formula above in the form we will need it for
our study of zeta functions.  However, from the toric geometry
point of view, it is more natural to use the continued fraction
expansion of $\ds\frac{q}{q-p}$ instead.  This formula will
be given in Section \ref{sub.proof12}.
\end{remark}

\begin{remark}
The formula for the Todd operator in Theorem \ref{thm.toddd} is
reminiscent of formulas in quantum cohomology, even though we know
of no conceptual explanation for this fact.
\end{remark}

\subsection{Dedekind sums}
\lbl{sub.dedekind}

Finally, we calculate the coefficients of the Todd power series of 
two-dimensional cones in terms of a  particular generalization $s_{i,j}$
(defined below) of  the classical Dedekind sums.
For an excellent review of the properties
of the classical Dedekind sum, see \cite{R-G}. Several
generalizations of the classical Dedekind sums were studied by T. Apostol, 
Carlitz, C. Meyer, D. Solomon,
and more recently and generally by U. Halbritter \cite{AV, Ha, Me, So}.
These papers investigate the sums ( and several other generalizations of
them, which we will not consider here)
given in  the following definition:
\begin{definition}
\lbl{def.dedsum}
For relatively prime integers $p,q$ (with $ q \neq 0$), and for
nonnegative integers $i,j$, we define 
the following {\em generalized Dedekind sum}:
\begin{equation}
\lbl{eq.dsum}
s_{i,j}(p,q)= \frac{1}{i! j!} 
\sum_{a=1}^{q}\hB_i(\lp \frac{a}{q} \rp) \hB_j(\lp -\frac{ap}{q} \rp)
\end{equation}
where for a real number $x$, we denote by
$\lp x \rp $ the (unique)  number such that $\lp x \rp \in (x + \BZ) \cap
(0,1]$. Here $\Ber_m$ denotes the $m^{th}$
{\em Bernoulli polynomial}, defined by the power series
$ \sum_{m=0}^{\infty} \Ber_m (x) \frac{t^m}{m!} =\frac{ t e^{tx}}{1 - e^t}$,
and $\hB_m$ denotes the restriction of the $m^{th}$ Bernoulli polynomial
$\Ber_m$ to $(0,1]$, with the boundary condition $\hB_m(1) \eqdef
\frac{1}{2} ( \Ber_m(1) + \Ber_m(0)) = B_m + \delta_{m,1}/2$, where
$B_m$ is the $m^{th}$ Bernoulli number, defined by $B_m=\Ber_m(0)$.
\footnote{There seem 
to be two conventions for denoting Bernoulli numbers; one, which we follow 
here, is used  often in number theory texts.
Due to the facts that $B_{2n+1}=0$ for $n\geq 1$, and  that
the sign of $B_{2n}$ alternates, the other 
convention, which is often used in intersection theory, 
defines the $n^{th}$ Bernoulli number to be $(-1)^{n+1} B_{2n}$.}
\end{definition} 
By definition, the sum $s_{1,1}(-p,q)$ coincides
with the classical Dedekind sum $s(p,q)$ (cf. \cite{R-G}.)
An  important property of generalized Dedekind sums is a reciprocity formula
which leads to an evaluation in terms of negative-regular continued 
fractions.  This reciprocity formula was most conveniently
 written by D. Solomon in terms of an additivity formula of a  generating 
power series 
\cite[Theorem 3.3]{So}. On the other hand, the Todd operator also satisfies 
an additivity property. 
We denote by $f_{i,j}(p,q) $ (for nonnegative integers $i,j$)
the coefficient of $x^i y^j$ of the power series 
$\ft_{\s_{(p,q)}}(x,y)$
({\em abbreviated} by $\ft_{(p,q)}(x,y)$).
We then have:

\begin{theorem}
\lbl{thm.f=s}
Let $p,q\in\BZ$ be relatively prime with $q \neq 0$. If $i,j>1$, then we have:
\begin{equation}
\lbl{eq.spp}
f_{i,j}(p,q)=q^{i+j-1} (-1)^i
 s_{i,j}(p,q).
\end{equation}
If $i=1$ or $j=1$, the above equation is true when the correction term
$$
q^{i+j-1}(-1)^{i+j}\frac{B_iB_j}{i!j!}
$$
is added to the right hand side.
\end{theorem}

\begin{corollary}
Fixing $r$, 
for all nonnegative integers $i,j$ and with the notation of \eqref{eq.bis}, 
$s_{i,j}(p,q)$ are polynomials in $b_i, 1/q$.
\end{corollary}

\begin{remark}
For $i=j=1$, the above theorem was obtained in 
\cite{P1}, and was a motivation for the results of the present paper.
\end{remark}

\subsection{Is toric geometry needed?}
\lbl{sub.needed}

Some  natural questions arise, at this point: 
\begin{itemize}
\item
Is toric geometry needed?
\item
How do the statement and proof of Theorem \ref{cor.zv} differ from the 
statement and proof of Zagier's formula  \cite{Za2}?
\end{itemize} 

With respect to the first question, Theorem \ref{cor.zv} 
(an evaluation of zeta functions associated to real quadratic fields)
is stated without
reference to toric geometry, and a close examination reveals that
its proof is based on an analytic Lemma \ref{prop.1} and on the two
dimensional analogue of the Euler-MacLaurin formula given by Proposition
\ref{prop.11}. In addition, two dimensional cones can be canonically
subdivided into nonsingular cones, so in a sense the two dimensional analogue
of the Euler-MacLaurin formula can be obtained by the classical 
(one dimensional) one, as is used by Zagier \cite{Za2}. 
Furthermore, number theory offers, for every $m \geq 2$, a canonical 
{\em Eisenstein cocycle}
of $\text{PGL}_m(\BQ)$ \cite{Sc1, Sc2} that expresses, among other things,
the generalized  Dedekind sums $s_{i,j}$ in terms of
negative-regular continued fraction expansions like the ones of Theorems
\ref{thm.toddd} and  \ref{thm.f=s}, and the values (at nonpositive integers)
of zeta functions of totally real fields (of degree $m$) in terms
of generalized Dedekind sums. See also \cite{So, St}.

On the other hand, toric geometry constructs for every simplicial cone
$\s$ in an $m$-dimension- al lattice $N$, a canonical Todd power series
$\ft_\s$ satisfying an $N$-additivity property (see Proposition \ref{prop.j1}
below). 
The coefficients of the Todd power series are generalized Dedekind sums,
and the power series itself is intimately related to the Euler-MacLaurin
summation formula. As a result, for $m=2$, we provide a toric geometry
explanation of Theorems \ref{cor.zv} and \ref{thm.f=s}. 

With respect to the second question, Zagier \cite{Za2}
obtained similar formulas for the values of the zeta function of a quadratic 
number field at negative integers. Zagier used additivity in the
lattice $M$, whereas
we use additivity in the dual lattice $N$. $N$-additivity 
implies  that these zeta values
are  a priori polynomials in $b_i$ and $1/q$; moreover these polynomials are
directly defined in terms of the Todd operator and the quadratic form $Q$.
Thus the present paper is an example of the rather subtle
 difference between $N$-additivity and $M$-additivity. 
This difference can be seen by contrasting the additivity
property of the zeta function $\zeta_{Q,\tau}$ with the additivity
property satisfied by Todd operator.
The modified Todd operator $\fs_\s$ (introduced in Section \ref{sub.general})
is $N$-additive
under subdivisions of the cone $\sigma$, i.e.,
it satisfies equation \eqref{eq.j1} below. On the other hand, the function
$\z_{Q, \cdot}(s): \C^2(M) \to \BC$ is $M$-additive, i.e., it satisfies:
\begin{equation}
\lbl{eq.Mad}
\z_{Q,\tau}(s) = \z_{Q, \tau_1}(s) +
\z_{Q, \tau_2}(s)
\end{equation}
whenever $\tau$ is subdivided into cones $\tau_1$ and $\tau_2$.
This $M$-additivity
of $\z_{Q, \cdot}(s)$ is due to the particular choice of the weight function
$wt$,  and is purely set-theoretic. In contrast, the $N$-additivity of
$\fs_{\s}$ is somewhat deeper, and has the advantage that only  cones of
 maximal dimension are involved.

\subsection{Plan of the proof}
\lbl{sub.plan}

In Section \ref{sec.thm.1}, we review well known properties of Todd
power series and prove Theorems \ref{thm.1} and \ref{thm.toddd}.
In Section \ref{sec.calc}, we review the relation between the zeta functions
that we consider and the zeta functions of real quadratic fields. We
give a detailed construction of zeta functions, and  prove Theorem 
\ref{cor.zv}. Finally, in Section \ref{sec.dedsum} we review properties
of generalized Dedekind sums and prove Theorem \ref{thm.f=s}.

\subsection{Acknowledgments}

We wish to thank M. Rosen for encouraging conversations
during the academic year 1995-96. We also thank B. Sturmfels
and M. Brion for their guidance. We especially wish to thank 
W. Fulton
for enlightening and encouraging conversations since our early years
of graduate studies.

\section{The Todd power series of a cone}
\lbl{sec.thm.1}
In this section we study properties of the Todd power series
$\ft_{\sigma}$ associated to a simplicial cone $\sigma$.
Section \ref{sub.general} provides an introduction and statements of 
several previously discovered 
formulas for the Todd power series.  In Section \ref{sub.2dim}
we make these formulas explicit for two-dimensional cones.  Section 
\ref{sub.proof12} contains a proof of the explicit continued fraction
formula for the Todd power series of a two-dimensional cone.  Finally,
in Section \ref{sub.thm.1} we prove Theorem \ref{thm.1} which links
Todd power series with the problem of evaluating zeta functions
an nonpositive integers.

\subsection{General properties of the Todd power series}
\lbl{sub.general}

Todd power series were studied in connection with the Todd class
of a simplicial toric variety in \cite{P3}.  Independently, they were 
introduced in \cite{BV2} in the study 
of the equivariant Todd class of a simplicial toric variety.  In
addition, these power series appear in Brion and Vergne's formula
for counting lattice points in a simple polytope \cite{BV}, which
is an extension of Khovanskii and Pukhlikov's formula \cite{KP}.  
In these remarkable formulas, the power series in
question are considered as differential operators which are applied
to the volume of a deformed polytope.  The result yields the number
of lattice points in the polytope, or more generally, the sum of any
polynomial function over the lattice points in the polytope.  Below
(Proposition \ref{prop.11}), we give a version of this formula expressing the
sum of certain functions over the lattice points contained in a simplicial
cone.

The Todd power series considered in the works cited above
are also closely related to the fundamental work of R. Morelli
on the Todd class of a toric variety.
\cite{Mo}.  A precise connection is given in \cite[Section 1.8]{P3}.

We now state some of the properties of the Todd power series of 
a simplicial cone.  Our purpose is twofold: we will need these properties
in our application to zeta functions, and we wish to unite the
approaches of the works cited above.  Here, we follow the notation of
\cite{P3}.

Let $\sigma=\lp \rho_1, \dots \rho_n\rp$ be an $n$-dimensional
simplicial cone in an $n$-dimensional lattice $N$.  The Todd 
power series $\ft_{\sigma}$ of $\sigma$ is a power series with
rational coefficients in
variables $x_1, \dots,
x_n$ be corresponding to the rays of $\sigma$.  These power series, when
evaluated at certain divisor classes, yield the Todd class of any
simplicial toric variety \cite[Theorem 1]{P3}.  

To state the properties of $\ft_{\sigma}$, it will be useful to consider
the following variant $\fs$ of the power series $\ft$ defined in \cite{P3}
by:
$$
\fs_{\sigma}(x_1,\dots,x_n)=\frac1{\text{mult}(\sigma) x_1\cdots x_n}
\ft_{\sigma}(x_1,\dots,x_n),
$$
which is a Laurent series in $x_1, \dots, x_n$.

The Todd power series $\ft_{\sigma}$ and $\fs_{\sigma}$ are
characterized by the following proposition \cite[Theorem 2]{P3},
which states that $\fs$ is {\it additive} under subdivisions (after
suitable coordinate changes), and
gives the value of $\fs$ on nonsingular cones.
An $n$-dimensional cone is called nonsingular if it is generated
by rays forming a basis of the lattice.  It is well-known that
any cone may be subdivided into nonsingular cones, and that such
a subdivision determines a resolution of singularities of
the corresponding toric variety (cf. \cite[Section 2.6]{Fu}.)

\begin{proposition}
\lbl{prop.j1}
If $\Gamma$ is a simplicial subdivision of $\sigma$ then
\begin{equation}
\lbl{eq.j1}
\fs_{\sigma}(X)=\sum_{\gamma\in\Gamma_{(n)}}
\fs_{\gamma}(\gamma^{-1}\sigma X).
\end{equation}
where the sum is taken over the $n$-dimensional cones of the subdivision
$\Gamma$.
Here $X$ denotes the column vector $(x_1,\dots,x_n)^t$ and
we have identified each cone ($\sigma$ and $\gamma$) with the $n$-by-$n$ matrix
whose columns are the coordinates of the rays of that cone.

For nonsingular cones,
$\sigma$, we have the following expression for $\fs_{\sigma}$:
$$
\fs_{\sigma}(x_1,\dots,x_n)=\prod_{i=1}^n\frac{1}{1-e^{-x_i}}.
$$
\end{proposition}

\begin{remark} It follows immediately that for any cone $\sigma$, 
$\fs_{\sigma}$ is a rational function of the $e^{x_i}$.
\end{remark}

The Laurent series $\fs_{\sigma}$ in $x_1, \dots , x_n$
may also be expressed as an exponential
sum over the lattice points in the cone, in the spirit of the important 
earlier work of M. Brion
\cite{Br}. 

\begin{proposition}
\lbl{prop.j2}
 Let $\check\sigma$ denote the dual cone in the
lattice $M=\text{Hom}(N,\BZ )$.  We then have
$$
\fs_{\sigma}(x_1,\dots,x_n)=\sum_{m\in\check\sigma\cap M}
e^{-(\lp m,\rho_1 \rp x_1 + \cdots + \lp m,\rho_n \rp x_n )},
$$
The equality is one of rational functions in the $e^{x_i}$.
\end{proposition} 

\begin{proof}  As noted in \cite{Br}, since $\lp m,\rho_i \rp \geq 0$
for $m\in\check\sigma\cap M$, the right hand side has a meaning
in the completion of $\BC[y_1,\dots,y_n]$  with respect to the 
ideal $(y_1,\dots,y_n)$,  where $y_i$ stands for
$e^{-x_i}$. 
While it is not obvious, the right hand side is an rational function of the
$e^{x_i}$.  See \cite[p.654]{Br}.  The left hand side
is a rational function of the $e^{x_i}$ by Proposition \ref{prop.j1}.  
By \cite[p.655]{Br}
and the second formula of  Proposition \ref{prop.j1}, these
two rational functions are equal on nonsingular cones.  As
any cone can be subdivided into nonsingular cones, it suffices
to verify that the right hand side satisfies the additivity formula
of Proposition \ref{prop.j1}. But this follows again from Brion's work. 
See the proposition of \cite[p.657]{Br}, for example.  Intuitively,
this additivity can be seen from the fact that a sum of exponentials
over a cone containing a straight line vanishes formally.
\end{proof}

The $\fs_{\sigma}$ also have an explicit expression in
terms of cyclotomic sums, due to Brion and Vergne.  Following
\cite{BV2}, we introduce the following notation.  Let $u_1,\dots,
 u_n$ denote the primitive generators of the dual cone $\check\sigma$.
Thus we have $\lp u_i,\rho_j \rp =0$ if $i\ne j$, and $\lp u_i,\rho_i \rp \in
\BZ $, but does not necessarily equal $1$.   Let $N_{\sigma}$ be the subgroup
of $N$ generated by the $\rho_i$, $i=1,\dots,n$, and let $G_{\sigma}=
N/N_{\sigma}$.  Then $G_{\sigma}$ is an abelian group of order
$\text{mult}( \sigma)$, which we denote by $q$.  
Define characters $a_i$ of $G_{\sigma}$ by 
$$
a_i(g)=e^{2\pi i\frac{\lp u_i, g \rp}{\lp u_i, \rho_i \rp}}.
$$

\begin{proposition}
\lbl{prop.j3}
The Laurent series $\fs_{\sigma}$ coincides
with Brion and Vergne's formula expressing their Todd differential operator.
Namely, we have
$$
\fs_{\sigma}(x_1,\dots,x_n)=\frac 1q \sum_{g\in G_{\sigma}}\prod_{i=1}^n
\frac 1{1-a_i(g)e^{-x_i}}.
$$
\end{proposition}

\begin{proof} By inspection, the right hand side is a rational function
in the variables $y_i=e^{-x_i}$ which takes the value $1$ when $y_i=0$
for all $i$. Such rational functions embed
into the completion of the ring $\BC[y_1,\dots,y_n]$ 
with respect to the ideal $(y_1, \dots , y_n)$.  To prove the proposition,
it is enough to show that the right hand side and the right hand side
of Proposition \ref{prop.j2} define the same element of this completion.  
To do so, we expand the right hand side, getting:
$$
\frac 1q \sum_{k_1,\dots,k_n\ge 0}
\sum_{g\in G_{\sigma}}
\prod_{i=1}^n
\biggl[a_i(g)e^{-x_i}\biggr]^{k_i},$$
which becomes
$$
\frac 1q\sum_{k_1,\dots,k_n\ge 0}
e^{-(k_1 x_1 +\cdots + k_n x_n)}
\sum_{g\in G_{\sigma}}
a_1^{k_1}\cdots a_n^{k_n}(g).
$$
This last sum over $G_{\sigma}$ is either $q$ or $0$, depending on whether
$a_1^{k_1}\cdots a_n^{k_n}$ is the trivial character of $G_{\sigma}$
or not.  Comparing the above with the right hand side of Proposition
2, we see that it suffices to show that $a_1^{k_1}\cdots a_n^{k_n}$ is the 
trivial character of $G_{\sigma}$ if and only if there exists $m\in
\check\sigma$ such that $\lp m,\rho_i \rp = k_i$ for all $i$. This
straightforward lattice calculation is omitted. 
\end{proof}

\begin{remark}
The sum in the proposition above appears in the important work
of Ricardo Diaz and Sinai Robins
\cite{DR}.  They use such sums to give an
explicit formula for the number 
of lattice points in a simple polytope.  Interestingly, their techniques,
which come from Fourier analysis, are seemingly unrelated to the toric
geometry discussed above.
\end{remark}

\subsection{Properties of Todd power series of two-dimensional cones}
\lbl{sub.2dim}

In this section, we state properties of the power series $\ft_{\sigma}$
for a two-dimensional cone $\sigma$. 

It is not hard to see that if $\sigma$ is
any two-dimensional cone, then there are relatively prime
integers $p,q$ such that
$\sigma$  is lattice-equivalent
to the cone
$$
\s_{(p,q)} \eqdef \lp (1,0), (p,q) \rp \subset \zz^2.
$$
Here $q$ is determined up to sign and $p$ is determined
modulo $q$.  Thus we may arrange to have $q>0$ and $0\leq p<q$.
This gives a complete classification of two-dimensional cones up to lattice
isomorphism.  In discussing Todd power series we will abbreviate
$\ft_{\s_{(p,q)}}$ by $\ft_{(p,q)}$.

The explicit cyclotomic formula of Proposition \ref{prop.j3} may be 
written as:

\begin{proposition}
\lbl{prop.2d.cycl}
The Todd power series of a two-dimensional cone is given by
$$
\ft_{(p,q)}(x,y)= 
\sum_{\w^q=1} \frac{xy}{(1- \w^{-p} e^{-x})(1 - \w e^{-y})}
$$
\end{proposition}

\begin{proof}
With the coordinates above and the notation of Proposition \ref{prop.j3}, 
we have $u_1=(q,-p)$, $u_2=(0,1)$, and $G_{\sigma}$ consists of
the lattice points $(0,k), k=0,\dots, q-1.$  The desired
equation now follows directly from Proposition \ref{prop.j3}.
\end{proof}

In the two-dimensional case, 
the additivity formula of Proposition \ref{prop.j1}
can be expressed as an explicit {\em reciprocity} law.
This, and a 
{\em periodicity} relation for $\fs$ are contained in the following theorem.

\begin{proposition}
\lbl{prop.2d.reci}
Let $p$ and $q$ be relatively prime positive integers.  Then
\begin{eqnarray}
\lbl{eq.reci}
\fs_{(p,q)}(x - \frac{p}{q} y, \frac{1}{q} y)+
\fs_{(q,p)}(y- \frac{q}{p}x , \frac{1}{p} x) & = &
\fs_{(0,1)}(x,y) \\
\lbl{eq.peri}
\fs_{(p+q,q)}(x,y) & = & \fs_{(p,q)}(x,y)
\end{eqnarray}
\end{proposition}

\begin{proof} The quadrant $\lp (1,0),(0,1) \rp$ may be subdivided
into cones $\gamma_1=\lp (1,0), (p,q) \rp$ and 
$\gamma_2=\lp (0,1), (p,q) \rp$.  The cone $\gamma_1$ is of type
$(p,q)$, and $\gamma_2$ is of type $(q,p)$.  Applying the additivity
formula (Proposition \ref{prop.j1}) to this subdivision yields the 
first equation above.

The second equation follows from the fact that the cones
$\lp (1,0), (p,q) \rp$ and
$\lp (1,0), (p+q,q) \rp$
are lattice isomorphic.
\end{proof}

Let  $\fs^{ev}$
(resp. $\fs^{odd}$) denote the part of $\fs$ of even (resp. odd) total
degree.

\begin{corollary}
\lbl{cor.oneonly}
Given a two-dimensional lattice $N$, the function 
$\fs: \C^2(N) \to \BQ\(x,y\)$ (where $\BQ\(x,y\)$ is the function field
of the power series ring $\BQ\[x,y\]$) is uniquely determined by
properties \eqref{eq.reci} and \eqref{eq.peri} and its {\em initial
condition} $\fs_{(0,1)}(x,y)= 1/((1 - e^{-x})(1 - e^{-y}))$.
Since equations \eqref{eq.reci} and \eqref{eq.peri} are homogenous with
respect to the degrees of $x$ and $y$, it follows that  $\fs^{ev}$
(resp. $\fs^{odd}$) satisfies \eqref{eq.reci} and \eqref{eq.peri}
with initial condition  $\fs_{(0,1)}^{ev}$ (resp. $\fs_{(0,1)}^{odd}$). 
\end{corollary}

\begin{remark}
\lbl{rem.oneo}
Corollary \ref{cor.oneonly} has a converse.
>From equations \eqref{eq.reci} and \eqref{eq.peri} it follows that
$\fs_{(0,1)}$ satisfies the following relations:
$$
\fs_{(0,1)}(x,y)  =  \fs_{(0,1)}(x-y,y) + \fs_{(0,1)}(y-x, x) \text{ and }
\fs_{(0,1)}(x,y)  =  \fs_{(0,1)}(y,x)
$$
Conversely, 
one can show that given any element  $g(x,y) \in \BQ\(x,y\)$
satisfying:
$$
g(x,y)  =  g(x-y,y) + g(y-x, x) \text{ and }
g(x,y)  =  g(y,x)
$$ 
there is a unique function $\mathfrak g :\C^2(N) \to \BQ\(x,y\)$
so that $\mathfrak g_{(0,1)} = g$.
\end{remark}

\subsection{Continued fraction expansion for the Todd power series}
\lbl{sub.proof12}

In this section, we prove Theorem \ref{thm.toddd}, which
expresses the coefficients in the Todd power series of
a two-dimensional cone in terms of continued fractions.
Before doing so, we first formulate an equivalent
version which is more natural from the point of view
of toric varieties.  The continued fraction expansion
in this second version of the formula corresponds directly
to a desingularization of the cone. 

Given relatively prime integers $p,q$ with $p>0$ and $0\leq 
p<q$, let
$ a_i,  \num_i,  \denom_i, L_i$  be defined in terms of the negative continued
fraction expansions:
\begin{equation*}
\frac{q}{q-p}  =  [a_1,\dots, a_{s-1}], \qquad
\frac{ \num_i}{ \denom_i}  \eqdef  [a_1,\dots,a_{i-1}], \qquad
 L_i  \eqdef   \num_i x+(q \denom_i+(p-q) \num_i)y.
\end{equation*}

We then have the following continued fraction expression for the degree $d$  
part $(\ft_{\conevar})_d$  of the Todd power series of a two-dimensional
cone $\conevar$.

\begin{theorem}
\lbl{thm.todd}
For $\conevar$ a cone of type $(p,q)$ as above, and for  $d\geq 
2$  an even integer, we have:
\begin{equation*}
(\ft_{\conevar})_{d}(x,y)=-qxy \sum_{i=1}^s
P_d(L_{i-1}, L_i) -\lambda_d qxy\sum_{i=1}^{s-1} a_i 
R_d( L_{i-1}, L_{i+1})
+\lambda_d(x L_1^{d-1}-y  L_{s-1}^{d-1})
\end{equation*}
If $d\geq 1$ is odd, then $
(\ft_{\conevar})_{d}(x,y)
=\frac12\lambda_{d-1}q^{d-1}xy(x^{d-2}+y^{d-2}). $
\end{theorem}

\begin{proof}
It will be convenient to 
choose coordinates so that $\conevar= \lp (0,-1), (q,q-p) \rp$
(which is easily seen to be lattice equivalent to the cone 
$\lp (1,0),(p,q) \rp$). We now subdivide $\conevar$ into
nonsingular cones. It is well-known that for two-dimensional
cones this can be done in a canonical way, and that the resulting
subdivision has an explicit expression in terms of continued fractions
\cite[Section 2.6]{Fu}.
In our coordinate system, the rays of this nonsingular subdivision of $\conevar$
are given by 
\begin{eqnarray*}
 \beta_0&=&(0,-1)\\
 \beta_1&=&(1,0)\\
 \beta_2&=&(a_1,1)\\
&\dots&\\
 \beta_s &=&(q,q-p).
\end{eqnarray*}
Thus we  have $ \beta_{i+1} +  \beta_{i-1} = a_i  \beta_i$, 
which implies that
 $ \beta_i=( \num_i,  \denom_i)$.

The cone $\conevar$ is subdivided into cones $\conevar_i= \lp  \beta_{i-1}, 
\beta_i \rp$,
$i=1,\dots s$.  The $N$-additivity formula of Proposition \ref{prop.2d.reci}
implies that:
$$
\fs_{\conevar}(x,y)=\sum_{i=1}^s \fs_{\conevar_i} (\conevar_i^{-1}\conevar (x,y)^t),
$$
where we again we have identified the $2$-dimensional cone $\conevar$ with 
the 2-by-2 matrix
whose columns are the primitive generators of $\conevar$.
One easily sees that this becomes
$$
\fs_{\conevar}(x,y)=\sum_{i=1}^s \fs_{\conevar_i} (-L_{i-1}, L_i)
$$
Rewriting the equation in terms of the $\ft_{\conevar}$ yields
$$
\ft_{\conevar}(x,y)=-qxy\sum_{i=1}^s \frac{\ft_{\conevar_i} 
(-L_{i-1}, L_i)}{L_{i-1}L_i}
$$
Since every $\conevar_i$ is nonsingular, we have
$$
\ft_{\conevar_i}(X,Y)=g(X)g(Y)
$$
where
$g(z)=z/(1-e^{-z}).$
Now  consider the degree $d$ part of the above. We will assume
$d>2$ is even and leave the other (easier) cases to the reader.
The degree $d$ part of $g(L_{i-1})g(L_{i})$ may be written as:
$$
-L_{i-1}L_iP_d(-L_{i-1},L_i) + 
\lambda_d (L_{i-1}^d+L_i^d).
$$
Summing the first term above yields  the first term in the
equation of the theorem.  So it suffices to examine the remaining
term:
$$
-\lambda_dqxy\sum_{i-1}^s \frac{L_{i-1}^d+L_i^d}{L_{i-1}L_i}.
$$
Using the relation $L_{i-1}+L_{i+1}=a_iL_i$, the sum above may be
rewritten as
$$
\sum_{i=1}^{s-1} a_i \frac{L_{i-1}^{d-1}+L_{i+1}^{d-1}}{L_{i-1}+L_{i+1}}
+\biggl(\frac{L_1^{d-1}}{L_0} + \frac{L_{s-1}^{d-1}}{L_s} \biggr).
$$
Keeping in mind that $L_0=-qy$ and $L_s=qx$, Theorem \ref{thm.todd}
follows easily.
\end{proof}

\begin{corollary}
\lbl{cor.use1}
For a two-dimensional cone $\s= \lp \rho_1, \rho_2 \rp$ of multiplicity $q$
we have:
\begin{equation}
\lbl{eq.evenf}
\ft_{\s}(h_1, h_2) - \frac{1}{2}\ft_{\lp \rho_1 \rp}( q h_1)h_2
- \frac{1}{2}\ft_{\lp \rho_2 \rp}( q h_2) h_1
= \ft_{\s}^{ev}(h_1, h_2) - \frac{q  h_1 h_2}{2}
\end{equation}
where $\ft_{\s}^{ev}$ is the even total degree part of the power series
$\ft_{\s}$.
\end{corollary}

\begin{proof}
It follows immediately from (and in fact is equivalent to)
 the formula for the odd part of the Todd power 
series $\ft_{\s}$ given by  Theorem \ref{thm.todd}.
\end{proof}

We now prove Theorem \ref{thm.toddd}.
Let $\sigma$ be a two-dimensional cone of type $(p,q)$
in a lattice $N$, and let $ a_i,  h_i,  k_i$ and $  X_i$
be as in Theorem \ref{thm.toddd}. The dual cone $\check{\sigma}$ in the 
dual lattice $M$ is easily seen to have type $(-p,q)$ and so we may
choose coordinates in $M$ so that
$\check{\sigma}= \lp (0,-1), (q,p) \rp$ in $\BZ^2$. 
Furthermore, the negative-regular continued fraction expansion of $q/p$
corresponds naturally to the desingularization of the dual cone 
$\check{\sigma}$.  Explicitly, the desingularization of $\check\sigma$
is given by the subdivision
\begin{eqnarray*}
 \rho_0&=(0,-1)\\
 \rho_1&=(1,0)\\
 \rho_2&=(b_1,1)\\
&\dots\\
 \rho_r &=(q,p).
\end{eqnarray*}
One has $ \rho_{i+1} +  \rho_{i-1} =b_i\rho_i$ and thus
 $ \rho_i=( h_i, k_i)$. Applying Theorem \ref{thm.todd}
to $\cs$ expresses $\ft_{(-p,q)}$ in terms of the $b_i$.  However
 $\ft_{(-p,q)}$ is related to $\ft_{(p,q)}$ via the
relation
$$\ft_{(p,q)}(x,y)= \ft_{(-p,q)}(-x,y) + \frac{qxy}{1-e^{-qy}}.$$
In this way, we obtain an expression for $\ft_{(p,q)}$ in terms of the $b_i$,
which concludes the proof of Theorem \ref{thm.toddd}.

\subsection{Zeta function values in terms of Todd power series}
\lbl{sub.thm.1}

In this section, we prove Theorem \ref{thm.1} which expresses
values of the zeta function of a two-dimensional cone in terms of
the Todd power series.
Three ingredients are involved in the proof of this theorem: an
asymptotic series formula (Proposition \ref{prop.1}), a polytope
summation formula (Proposition \ref{prop.111}) and a cone summation
formula (Proposition \ref{prop.11}).
We begin with the first ingredient:
 
\begin{lemma}\cite[Proposition 2]{Za2}
\lbl{prop.1}
Let $\phi(s)= \sum_{\l > 0} a_\l \l^{-s}$ be a Dirichlet series
where $\{ \l \}$ is a sequence of positive real numbers converging to infinity.
Let $E(t)= \sum_{\l > 0} a_\l e^{- \l  t}$ be the corresponding 
exponential series. Assume that $E(t)$ has the 
following asymptotic 
expansion as $t \to 0$:
\begin{equation}
E(t) \sim  \sum_{n=-1}^{\infty}
c_{n} t^n
\end{equation}
Then it follows that 
\begin{itemize}
\item
  $\phi(s)$ can be extended to a meromorphic function on $\BC$.
\item
  $\phi(s)$ has a simple pole at $s=1$, and no other poles.
\item
  The values of $\phi$ at nonpositive  integers are given by:
$\phi(-n)= (-1)^{n} n! c_{n}$.
\end{itemize}
\end{lemma}

We now present our second ingredient, a polytope summation formula.  This
proposition is a variant of the lattice point formula of \cite{BV2}.
A polytope of dimension $n$ is called {\em simple} if each of its vertices lies
on exactly  $n$ facets ($n-1$-dimensional faces) of the polytope.
\begin{proposition}
\lbl{prop.111}
Let $N$ be an $n$-dimensional lattice,  $P$ a simple
 lattice polytope
in $M$ and $\Sigma$ its associated fan in $N$.
For every analytic function
$\phi: M_\BR \to \BR$, we have the following asymptotic expansion
as $t \to 0$:
\begin{equation}
\sum_{a \in P \cap M} \phi(t a) \sim
 \ft_{\Sigma}\left(\frac{\pa}{\pa h}\right) \di  \int_{P(h)} \phi(t u) du
\end{equation}
where $\ft_\Sigma$ is the Todd power series of \cite[Definition 10]{BV2}.
\end{proposition}

\begin{proof}
First of all, the meaning of the right hand side is as follows:
we consider the 
degree $k$ Taylor expansion $\phi=\phi_k + R_k$
of $\phi$, where $\phi_k$ is a polynomial in $M_\BR$ of degree $k$ and
$R_k$ is the remainder satisfying $\lim_{a\to 0} |a|^k R_k(a)=0$.
It follows that $\int_{P(h)} P_k(tu) du$ is a polynomial in $t$ and $h$ of
degree $k$ (with respect to $t$)  and that
$ \int_{P(h)} R_k(tu) du=o(t^k)$ at $t=0$ 
(with the notation that $f(t)=o(t^k)$ if and only if $\lim_{t\to 0}
f(t) t^{-k}=0$).
Thus, 
$$
\ft_\Sigma\left(\frac{\pa}{\pa h}\right)
\di \int_{P(h)} \phi(t u) du =
\ft_\Sigma\left(\frac{\pa}{\pa h}\right)
\di \int_{P(h)} \phi_k(t u) du + o(t^k).
$$
On the other hand, 
\begin{eqnarray*}
\sum_{a \in P \cap M} \phi(t a) & = &
\sum_{a \in P \cap M} \phi_k(t a) +
\sum_{a \in P \cap M} R_k(t a) \\
& = &   \sum_{a \in P \cap M} \phi_k(t a) + o(t^k).
\end{eqnarray*}
Brion-Vergne \cite[theorem 11]{BV2} prove that for every polynomial
function (such as $\phi_k$) on $M_\BR$ we have:
$$
\sum_{a \in P \cap M} \phi_k(t a) =
 \ft_{\Sigma}\left(\frac{\pa}{\pa h}\right) \di  \int_{P(h)} \phi_k(t u) du,
$$
which concludes the proof.
\end{proof}

We call a function $\phi: \BR^n \to \BR$ {\em rapidly decreasing}
if it is analytic and for every constant coefficients differential operator
$D$, and every subset $I$ of $[n]=\{1,\dots, n\}$, the restriction
$D(\phi)|_I$ obtained by setting $x_i=0$ for $i \not\in I$ is in 
$L^1(\BR_+^{I})$. Examples of rapidly decreasing functions
can be obtained by setting $\phi=P \exp(Q)$ where $P$ is a polynomial on 
$\BR^n$ and $Q: \BR^n \to \BR$ is totally positive, i.e., its
 restriction to $\BR_{+}^{I}$ takes positive values for every subset $I$ 
of $[n]$.

\begin{proposition}
\lbl{prop.11}
Let $N$ be an $n$-dimensional lattice.
For every $\s \in \C^n(N)$, and every rapidly decreasing analytic function 
$\phi: M_\BR \to \BR$,
we have the following asymptotic expansion as $t\to 0$:
\begin{equation}
\lbl{eq.p3}
\sum_{a \in \cs \cap M} \phi(t a) \sim
 \ft_{\s}\left(\frac{\pa}{\pa h}\right) \di  \int_{\cs(h)} \phi(t u) du
\end{equation}
\end{proposition}

\begin{proof}
First of all, the right hand side of the above equation has the
following meaning: consider the decomposition $\ft_\s=\sum_k \ft_{\s,k}$
of the power series $\ft_\s$, where $\ft_{\s,k}$ is a homogenous polynomial
of degree $k$. A change of variables $v=tu$ implies that
$$
\ft_{\s,k}\left(\frac{\pa}{\pa h}\right) \di  \int_{\cs(h)} \phi(t u) du
=
\ft_{\s,k}\left(\frac{\pa}{\pa h}\right) \di  \int_{\cs(th)} \phi(v) dv/t^n
=
t^{k-n}
\ft_{\s,k}\left(\frac{\pa}{\pa h}\right) \di  \int_{\cs(h)} \phi(v) dv
$$
is a multiple of $ t^{k-n}$, thus the right hand side is defined to be
the Laurent power series in $t$ given by 
$$
\sum_{k=0}^\infty t^{k-n} 
\ft_{\s,k}\left(\frac{\pa}{\pa h}\right) \di  \int_{\cs(h)} \phi(v) dv.
$$
\noindent
For the proof of the proposition, truncate in some way the cone $\cs$
in $M$ to obtain a simple convex polytope $P$, with associated fan $\Sigma$.
Since $\cup_{r > 0} r P=\cs$, using the convergence properties of $\phi$, we
obtain as $r \to \infty$:
\begin{eqnarray*}
\sum_{a \in \cs \cap M} \phi(t a) & = & \lim_r
\sum_{a \in r P \cap M} \phi(t a) \\ & \sim &
\lim_{r} 
 \ft_{\Sigma}\left(\frac{\pa}{\pa h}\right) \di  \int_{(rP)(h)} \phi(t u) du
\\ 
& = &
\lim_{r} 
 \ft_{\Sigma}\left(\frac{\pa}{\pa h}\right) \di  \int_{(rtP)(th)} \phi(v) 
dv/t^n \\
& = &
\lim_{r}
\sum_{k} \ft_{\Sigma,k}t^{k-n}\left(\frac{\pa}{\pa h}\right) \di
\int_{(rtP)(h)} \phi(v) dv \\
& = & 
\sum_{k} \lim_r
\ft_{\Sigma,k}t^{k-n}\left(\frac{\pa}{\pa h}\right) \di
\int_{(rtP)(h)} \phi(v) dv  \\
& = &  \sum_{k} 
\ft_{\s,k}t^{k-n}\left(\frac{\pa}{\pa h}\right) \di
\int_{P(h)} \phi(v) dv,
\end{eqnarray*}
which concludes the proof.
\end{proof}

In case of a two-dimensional lattice $N$ and a rapidly decreasing
function $\phi: M_\BR \to \BR$, using  the weight function $wt$
of equation \eqref{eq.we} and inclusion-exclusion, we obtain the following
 
\begin{corollary}
\lbl{cor.weights}
For a two-dimensional cone $\s  = \lp \rho_1, \rho_2 \rp$ of multiplicity
$q$ in $N$, we have the asymptotic expansion as $t \to 0$:
\begin{equation}
\lbl{eq.p3.weight}
\sum_{a \in \cs \cap M} wt(\cs, a) \phi(ta) \sim
\left\{ \ft_{\s}^{ev}\left(\frac{\p}{\p h_1},\frac{\p}{\p h_2}\right)  
-\frac{q}{2} \frac{\p^2}{\p h_1 \p h_2} \right\}
\di  \int_{\cs(h)} \phi(tu)du,
\end{equation}
where the right hand  side lies in the formal power series ring
$t^{-2} \BR \[t^2 \]$.
\end{corollary}

\begin{proof}[Proof of Theorem \ref{thm.1}]
Recall that $\tau$ is a cone in $M$, $\sigma$ is its dual in $N$ and
$Q$ is homogenous quadratic, totally positive on $\tau$; thus $e^{-Q}$
is rapidly decreasing on $M_\BR$. 
Corollary \ref{cor.weights} implies that the generating function
\begin{equation}
\lbl{eq.gf}
Z_{Q,\tau}(t)=
\sum_{ a \in \tau \cap M} wt(\tau, a)e^{-t Q(a)} 
=\sum_{ a \in \tau \cap M} wt(\tau, a)e^{-Q(t^{1/2}a)},
\end{equation}
satisfies the hypothesis
of Lemma \ref{prop.1},  which in turn yields Theorem \ref{thm.1}.
\end{proof}

\begin{remark}
\lbl{rem.notice}
Notice that the above proof of Theorem \ref{thm.1} used crucially
the fact that $B_1 = -\frac{1}{2}$ and the definition of the weight
function $wt$. If 
we had  weighted the sum defining the zeta function
in any other way, the
resulting variation of Theorem \ref{thm.1} would not hold.
\end{remark}

\section{Zeta functions of number fields}
\lbl{sec.calc}

\subsection{Parametrizing triples $(M_n,Q_b,\tau_b)$}
\lbl{sub.real}

In this section, we review in detail the construction 
of the triple $(M_b,Q_b,\tau_b)$ mentioned in Section \ref{sub.history}
 that is associated to
a sequence $b=(b_0, \dots, b_{r-1})$ of integers greater than $1$ and 
not all equal to $2$. Our notation and construction is borrowed from
\cite{Za2}. 
Given a sequence $b$ as above, we extend it 
to a sequence of integers parametrized by the integers by defining
$b_k = b_{k \bmod r}$. Furthermore, for an integer $k$, we 
define 
$$w_k= \[ b_k, \dots, b_{k+ r-1} \]
= b_k - \cfrac{1}{b_{k+1} - \cfrac{1}{b_{k+2} - \cdots }} 
$$
where $\[b_k, \dots, b_{k+ r-1} \]$ denotes the infinite continued fraction
with period $r$. Note that $w_k = w_{k+r}$ for all integers $k$, and that,
by definition, 
$$w_0= b_0 - \cfrac{1}{b_1 - \cfrac{1}{ \cdots b_{r-1} -\cfrac{1}{w_0}}} 
$$
from which it follows that $w_0$ satisfies  a  quadratic equation
$ A_b w^2 + B_b w + C_b = 0$ where $A_b,B_b,C_b$ are (for a fixed
$r$) polynomials in $b_i$
with integer coefficients.\footnote{not to be confused with
$A_0,\dots, A_{r}$ introduced below,
despite the rather poor choice of notation.}
 Let $D_b= B_b^2 - 4 A_b C_b$ be the discriminant.
Since $b_i > 1$, it follows that $D_b  > 0$, and that
 the roots of the quadratic equation are
$w_0=( -B_b + \sqrt{D_b})/(2 A_b)$ and $w_0 '=( -B_b - \sqrt{D_b})/(2 A_b)$.
We thus have that $w_0 > 1 > w_0 '$. We define the complex
numbers $A_0=1$, $A_{k-1}= A_k w_k$ for $k \in \BZ$. Since 
$w_k = b_k - \frac{1}{w_{k+1}}$ it follows that 
\begin{equation}
\lbl{eq.Ar}
A_{k-1} + A_{k+1}
= b_k A_k
\end{equation}
 for all integers $k$.  Let us define $M_b= \BZ w_0 + \BZ $
to be the rank two lattice 
 with the oriented basis $\{w_0, 1 \}$. Note that $M_b$ is a rank two
lattice in the  
real quadratic field $K_b = \BQ(\sqrt{D_b})$. Due to the fact that $A_{-1}
= w_0$ and the recursion relation \eqref{eq.Ar}, it follows that
$M_b= \BZ A_k + \BZ A_{k+1}$ for all integers $k$, and thus
$\epsilon M_b = M_b$ where $\epsilon = w_0 \dots w_{r-1}$.
Thus  $U_b= \{ \epsilon^m | m \in \BZ \} $ acts on $M_b$ by multiplication,
and the zeta function of (the inverse of) the narrow ideal class of
$M_b$ satisfies
\begin{equation*}
\z(M_b, s)= N(M_b)^s \sum_{\substack{a \in M_b/U_b \\
                       a \gg 0 }} \frac{1}{Q_b(a)^s},
\end{equation*}
where in the above summation we exclude $0$, and $N(M_b)$ is a constant,
defined in \cite{Za2}, which is not so crucial for our purposes. 
We are therefore interested to describe the quotient $M_b/U_b$.
The action of $U_b$ on $M_b$ implies  that $\ep A_k = A_{k-r}$.  
The matrix of multiplication by $\ep$ on $M_b$, in terms of the basis
$\{ \beta_1, \beta_2 \} \eqdef \{ A_{-1}, A_0 \}$ can be calculated as follows.
Recall first that in terms of the above basis of $M_b$,
we have that:
\begin{equation}
\lbl{eq.aaa}
A_k = - p_k A_{-1} + q_k A_0 \text{ where }
\frac{p_k}{q_k}= [ b_0, \dots, b_{k-1} ]
\end{equation}
which implies that
$$
\ep^{-1} \tByo {\b_1} {\b_2} = \tByo {A_{r-1}} {A_r} =
\tByt {-p_{r-1}} {- p_r} {q_{r-1}} {q_r} \tByo {\b_1} {\b_2} ,  
$$ 
from which it follows that the matrix of $\ep$ is given by
\begin{equation}
\lbl{eq.ea}
\tByt {q_r} {p_r} {-q_{r-1}} {-p_{r-1}} 
= \tByt
{ b_0 q - p} {q} {- b_0 p' + (p p' -1)/q} {- p'}
= \tByt {b_0} {-1} {1} {0} \dots  
\tByt {b_{r-1}} {-1} {1} {0}
\end{equation}
where  
\begin{equation}
\lbl{eq.pq}
\frac{q}{p} \eqdef [b_1, \dots, b_{r-1}], 
\end{equation}
(with the understanding that $q=1, p=0$ if $r=1$), 
$ p'=\text{numerator}
[b_1, \dots, b_{r-2}] $, which is a multiplicative inverse of
$p$ mod $q$ (with the understanding that 
$p'=0$ (resp. $1$) if $r=1$ (resp. $2$)).

\begin{lemma}
\lbl{claim.b}
The orbit space $M_b/U_b$ has as a fundamental domain a semiopen
two-dimensional simplicial cone whose closure is the cone
$\lp A_0,  A_r \rp$ in $K_b$. Furthermore, 
the (closed) cone $\tau_b=\lp A_0,  A_r \rp$ in $M_b$ is of type $(-p,q)$
and the dual cone $\sigma_b$
in $N_b$ is of type $(p,q)$ where $ 0 \leq p < q$ as in equation 
\eqref{eq.pq} above.
\end{lemma}

\begin{proof}
>From Figure \ref{Aray}, it follows that a fundamental domain of $M_b/U_b$ 
is a semiopen two-dimensional cone whose closure is 
the cone $\tau_b =\lp A_r,  A_0 \rp$ in $M_b$. 
The cone $\lp A_r,  A_0 \rp$
is canonically subdivided into nonsingular cones $\lp A_{i+1},  A_i \rp$
(for $i=0, \dots, r-1$), and using equation \eqref{eq.Ar}, it follows
that $\tau_b$ is of type $(c_1, q)$ where $0 \leq c_1 < q$ and 
$ \frac{q}{q - c_1} = [b_1, \dots, b_{r-1}]$.
Therefore the dual cone $\s_b$ in $N$ is of type $(p,q)$ where
$p= q-c_1$. Thus $(p,q)$ satisfies  equation \eqref{eq.pq}.
\end{proof}

\begin{figure}[htpb]
$$ \printname{Aray}
	\setlength{\unitlength}{0.03\standardunitlength}
	\begin{array}{c}  \hspace{-1.7mm}
        	\raisebox{-8pt}{\begingroup\makeatletter\ifx\SetFigFont\undefined
\def\x#1#2#3#4#5#6#7\relax{\def\x{#1#2#3#4#5#6}}%
\expandafter\x\fmtname xxxxxx\relax \def\y{splain}%
\ifx\x\y   
\gdef\SetFigFont#1#2#3{%
  \ifnum #1<17\tiny\else \ifnum #1<20\small\else
  \ifnum #1<24\normalsize\else \ifnum #1<29\large\else
  \ifnum #1<34\Large\else \ifnum #1<41\LARGE\else
     \huge\fi\fi\fi\fi\fi\fi
  \csname #3\endcsname}%
\else
\gdef\SetFigFont#1#2#3{\begingroup
  \count@#1\relax \ifnum 25<\count@\count@25\fi
  \def\x{\endgroup\@setsize\SetFigFont{#2pt}}%
  \expandafter\x
    \csname \romannumeral\the\count@ pt\expandafter\endcsname
    \csname @\romannumeral\the\count@ pt\endcsname
  \csname #3\endcsname}%
\fi
\fi\endgroup
\begin{picture}(2544,2133)(0,-10)
\thicklines
\path(12,12)(1512,612)
\path(1411.725,539.579)(1512.000,612.000)(1389.441,595.287)
\path(12,12)(1812,1212)
\path(1728.795,1120.474)(1812.000,1212.000)(1695.513,1170.397)
\path(12,12)(1512,1512)
\path(1448.360,1405.934)(1512.000,1512.000)(1405.934,1448.360)
\path(12,12)(612,1812)
\path(602.513,1688.671)(612.000,1812.000)(545.592,1707.645)
\path(12,12)(1812,312)
\path(1698.565,262.680)(1812.000,312.000)(1688.701,321.864)
\put(462,1962){\makebox(0,0)[lb]{$A_r$}}
\put(1737,612){\makebox(0,0)[lb]{$A_0= 1 = \b_2$}}
\put(2037,1137){\makebox(0,0)[lb]{$A_1$}}
\put(1662,1587){\makebox(0,0)[lb]{$A_2$}}
\put(612,1287){\makebox(0,0)[lb]{$\cdots$}}
\put(1962,162){\makebox(0,0)[lb]{$A_{-1}=w_0 = \b_1$}}
\end{picture} }
        	\hspace{-1.9mm}
	\end{array}
 $$
\caption{A fundamental domain of $M_b/U_b$.}\lbl{Aray}
\end{figure}
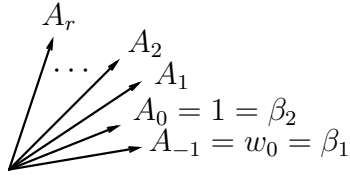

Since 
$M_b \subseteq \BQ(\sqrt{D_b})$ is a rank two lattice of a  real quadratic 
field,  we can define 
a homogenous function $Q_b: M_{\BR} \to \BR$ by 
$$Q_b(x w_0 + y)= C_b x^2  - B_b xy + A_b y^2
$$
This agrees with the (normalized) norm function of the real quadratic field 
$\BQ(\sqrt{D_b})$ considered by Zagier \cite[p 138]{Za2}.  
We close this section with the following:

\begin{proposition}
\lbl{prop.qij}
Let $l,m\in\BZ$ with $l\ne m$. Then the coefficients of $x^2, xy$ and $y^2$
of the quadratic form $Q_b(x A_l + y A_m)$ are  given by
polynomials in $b_i$ with integer coefficients.
\end{proposition}

\begin{proof}
We first claim that in terms of the basis $\{ A_r, A_0 \}$ we have:
\begin{equation}
\lbl{eq.change}
Q_b(x A_r + y A_0)= q Q_{\th}(x,y)
\end{equation}
where $Q_{\th}(x,y)= x^2 + \th xy + y^2$ and $\th= b_0 q - p - p'$.

The above claim follows from the fact that $A_b, B_b, C_b$ can be calculated
(using their definition) in terms of the $b_i$ as follows:
\begin{equation*}
A_b   =    q, \text{ \ \ \ }
B_b  =  - b_0 q + p - p', \text{ and }
  C_b   =   q_{r-1}= b_0 p' +(1-p p')/q,
\end{equation*}
together with a change of variables formula from $\{ A_{-1}, A_0 \}$
to  $\{ A_{r}, A_0 \}$ given by equation \eqref{eq.aaa}.

We then claim that if $l,m$ are distinct integers, then
\begin{equation}
\lbl{eq.qlm}
Q_b(x A_l + y A_m)= \frac{1}{q} Q_\Theta(H_l x + H_m y, K_l x + K_m y)
\end{equation}
where $\ds{ \qquad
\frac{ h_i}{ k_i} = [ b_1, \dots, b_{i-1}], \qquad
H_i x + K_i y  \eqdef   X_i \eqdef -  h_i x + (q  k_i -p  h_i)y. 
}$

This follows from equation \eqref{eq.change} using the  change of variables
given in equation \eqref{eq.aaa}. 

We now claim that the coefficients of $x^2, xy $ and $y^2$ in equation
\eqref{eq.qlm} are polynomials in the $b_i$ with integer coefficients.

Indeed, observe that $ h_i ,  k_i$ are polynomials in $b_i$
with integer coefficients. 
Since $H_i=- h_i, K_i = q  k_i -p  h_i = -p  h_i \bmod q$, 
it follows that:
\begin{eqnarray*}
Q_\Theta(H_l x + H_m y, K_l x + K_m y) & = &
( h_l x +  h_m y)^2 + \Theta ( h_l x +  h_m y)(p  h_l x + p
 h_m y) \\
& &
+ (p  h_l x + p  h_m y)^2 \bmod q \\
& = & ( h_l x +  h_m y)^2(1 + \Theta p + p^2) \bmod q. 
\end{eqnarray*}
Since $\Theta = q b_0 - p - p'$, and  $p p' \equiv 1  \bmod q$, 
 it follows that
$1 + \Theta p + p^2 = 0 \bmod q$, which finishes the proof of the proposition.
\end{proof}

\begin{remark}
\lbl{rem.compare}
The notation of Theorem \ref{thm.toddd} (and its proof,
in Section \ref{sub.proof12})
 corresponds to the
notation of this section. Indeed, the cone $\s_b$ in $N$ is of type
$(p,q)$ where $p,q$ are given in terms of equation \eqref{eq.pq}.
The rays $ \rho_i$ (for $i=1, \dots r-1$) that desingularize the
cone $\tau_b$ in $M$ are given by  $ \rho_i = A_i$. In addition, the
$ h_i,  k_i,  X_i$ that appear in Theorem \ref{thm.toddd}
match the notation  Proposition \ref{prop.qij}. 
\end{remark}

\subsection{Proof of Theorem \ref{cor.zv}}
\lbl{sub.explicit}

\begin{proof}
The main idea is to use Theorem \ref{thm.1} which calculates the zeta
values in terms of the Todd operator of $\s_b$, and Theorem \ref{thm.toddd}
which expresses the Todd operator in terms of the $b_i$.
The expression that we obtain for the zeta values differ from the one of
equation \eqref{eq.zva} by an error term, which vanishes identically,
as one can show by an explicit calculation.

Now, for the details, 
we follow the notation of Section \ref{sub.real}.
We begin by calculating 
the integral $\int_{\tau_b(h_1, h_2)} e^{-Q(u) } du$. 


Using the parametrization  $\BR^2 \to M_{\BR}$ given by:
$(w_1,w_2) \to w_1A_r+w_2A_0$, it follows that
the  preimage of $\tau_b(x, y)$ in $\BR^2$ is given
by $\{ (w_1, w_2) | w_1 \geq -x/q, w_2 \geq -y/q \}$.
Thus we have:  
\begin{eqnarray*}
 \int_{\tau_b(x, y)} e^{-Q_b(u) } du & = &
q \int_{w_1= -x/q}^{\infty} \int_{w_2= -y/q}^{\infty}
e^{-Q_b(q w_1, -p w_1 + w_2)} d w_2 d w_1
\end{eqnarray*}
Differentiating, 
we get
\begin{eqnarray*}
\left(\frac{\p}{\p x}\right) \left(\frac{\p}{\p y}\right)
\di \int_{\tau_b(x, y)} e^{-Q_b(u) } du & = 
\frac{1}{q} e^{-Q_b(- x A_r/q - y A_0/q)}
= \frac{1}{q} e^{-1/q Q_{\th}(x, y)}
\end{eqnarray*}
Using the following elementary identity:
$$
\left(\alpha \pd {\bar x} + \beta \pd {\bar y} \right)^i
\left(\gamma \pd {\bar x} + \delta \pd {\bar y} \right)^j \Big|_{\substack{
 \bar x = \alpha a + \gamma b \\
 \bar y = \beta a + \delta b }} f( \bar x, \bar y)=
\left(\pd x \right)^i \left(\pd y \right)^j 
\Big|_{\substack {x=a \\ y=b}} f(\alpha x +
\gamma y,  \beta x + \delta y)
$$
(and temporarily abbreviating $\pd x$ by $x$) we obtain that
\begin{eqnarray*}
\left(q xy  X_l^a  X_m^b \right) \di \int_{\tau_b(x, y)} e^{-Q_b(u) } du
& = &
 X_l^a  X_m^b \di e^{-1/q Q_{\th}(x, y)} \\
& = & \left(\pd x\right)^a\left(\pd y\right)^b 
\di e^{-1/q Q_\Theta(H_l x + H_m y, K_l x +K_m y)} \\
& = & \left(\pd x\right)^a\left(\pd y\right)^b 
\di e^{-Q_b(x A_l + y A_m)}, \quad\text{ as well as} 
\end{eqnarray*}
\begin{multline*}
 (x  X_1^{2n+1} + y  X_{r-1}^{2n+1}) \di  
\int_{\tau_b(x, y)} e^{-Q_b(u) } du = \\
   \frac{1}{q^{n+1}}
\int_{0}^{\infty} \left( \left(\pd x + p \pd y \right)^{2n+1} +
\left(\pd x + p' \pd y\right)^{2n+1} \right) \big|_{x=0} e^{-Q_{\th}(x,y)} dy
\end{multline*}

The above, together with Theorem \ref{thm.toddd}, implies that:
\begin{eqnarray*}
 \z_{Q_{b}, \tau_b}(-n) & = &
 (-1)^n n! \left\{ 
\sum_{i=1}^r P_{2n+2} \left(\pd x , \pd y \right) 
\di e^{- Q_b( x A_{i-1} + y A_i)}
\right. \\
& & \left.
+ \l_{2n+2} \sum_{i=1}^{r} b_i R_{2n+2}\left(\pd x , \pd y \right) \di   
e^{- Q_b( x A_{i-1} + y A_{i+1})} + E_{2n+2}(b)  \right\} 
\end{eqnarray*}
where 
\begin{eqnarray*}
E_{2n+2}(b)  & = & \l_{2n+2}
 \left\{ -b_0 R_{2n+2}\left(\pd x , \pd y \right) 
e^{-1/q(L x^2 + M xy + N y^2)} \right. \\ & & \left. -\frac{1}{q^{n+1}}
\int_{0}^{\infty} \left( \left(\pd x + p \pd y\right)^{2n+1} +
\left(\pd x + p' \pd y\right)^{2n+1} \right) 
\big|_{x=0} e^{-Q_{\th}(x,y)} dx \right\}
\end{eqnarray*}
and 
$$
L = q p' b_0 + 1 - p p' , \qquad
M = q b_0 \th_b +2(pp' -1), \qquad
N = q p b_0 + 1 - p p'.
$$
\newline
Fixing $r$, the length of the sequence $(b_0, \dots, b_{r-1})$, 
Lemma \ref{lem.zag} below can be used to express $E_{2n+2}(b)$ 
as a  polynomial 
in $1/q, p, p', b_0$. An explicit but lengthy calculation implies that
$E_{2n+2}(b)=0$ for any $b$. We have thus succeeded in expressing
$\z_{Q_{b}, \s_b}(-n)$
as a polynomial in $b_i$ with rational coefficients, 
symmetric under cyclic permutation, (since
$Q_b(x A_{l} + y A_m)$ is symmetric under cyclic permutation).  
This concludes the
proof of Theorem  \ref{cor.zv}.
Equation \eqref{eq.z0} follows easily using $\l_1= 1/2, \l_2=1/12$.
\end{proof}

\begin{lemma}\cite{Za2}
\lbl{lem.zag}
For every $i,j \geq 0$
with $i+j$ even, we have:
\begin{equation*}
\left(\frac{\p}{\p x}\right)^i \left(\frac{\p}{\p y}\right)^j 
\di e^{-(a x^2 + b xy + c y^2)}=
i! j! (-1)^{(i+j)/2} 
a^{i/2} c^{j/2} \sum_{i_2=0; i_2 \equiv i \bmod 2}^{min \{i,j \}}
\frac{ (a^{-1/2} b c^{-1/2})^{i_2}}{((i-i_2)/2)! i_2 ! ((j -i_2)/2)!}
\end{equation*}
Furthermore, for every $n \geq 1$
we have:
\begin{equation*}
\int_{y=0}^{\infty}
\left(\frac{\p}{\p x}\right)^{2n-1} \di
e^{-(a x^2 + b xy + c y^2)}  d y
=
- \frac{(2n-1)!}{2 c^{n}} 
\sum_{r=0}^{n-1} (-1)^r
\frac{(n-1-r)!}{r! (2n-1-2r)!} a^r b^{2n-1-2r} c^r
\end{equation*}
\end{lemma}

\section{Dedekind sums in terms of Todd power series}
\lbl{sec.dedsum}

In this section we give a proof of Theorem \ref{thm.f=s}.
To do this we will use Proposition \ref{prop.j2}, which
can be used to express the coefficients $f_{i,j}$ as
a sum of rational numbers, which turn out to equal products of
certain values of the Benoulli polynomials.  (Note, in contrast,
that Proposition \ref{prop.j3} expresses this same number in
terms of roots of unity instead of rational numbers.)

Let $p$ and $q$ be as in the statement of Theorem \ref{thm.f=s}, and 
let $\sigma$ be the cone $\lp (1,0),(p,q) \rp$ in $\zz^2$. We let
$\rho_1=(1,0)$, and $\rho_2=(p,q)$ denote the generators of this cone.
The left hand side of Theorem \ref{thm.f=s}, $f_{i,j}(p,q)$, equals the 
coefficient
of $x^iy^j$ in the power series $\ft_{\sigma}(x,y)$.  We now compute this 
power series using  Proposition \ref{prop.j2}.  This proposition
 contains an expansion for the power 
series $\fs_{\sigma}$, which is equivalent to the following expansion of 
$\ft_{\sigma}$:

$$
\ft_{\sigma}(x,y)=qxy\sum_{m\in\check\sigma\cap M}
e^{-(\lp m,\rho_1 \rp x + \lp m,\rho_2 \rp y )}.
$$

Every point of $\check\sigma\cap M$ can be written uniquely as a nonnegative
integral combination of the generators $u_1=(q,-p)$ and $u_2=(0,1)$ of 
$\check\sigma$, plus a lattice point in the semiopen parallelepiped
$$
P=\{cu_1+du_2 |  c,d\in [0,1) \}.
$$
Using $\lp u_1,\rho_1 \rp =\lp u_2,\rho_2 \rp = q$, it follows that

$$
\ft_{\sigma}(x,y)=qxy \frac1{1-e^{-qx}}\frac1{1-e^{-qy}}
\sum_{m\in P\cap M}
e^{-(\lp m,\rho_1 \rp x + \lp m,\rho_2 \rp y )}.
$$

One finds also that
$$P=\biggl\{(k,\{\frac{pk}{q}\}-\frac{pk}{q}) : k=0,\dots,q-1
\biggr\},
$$
where $\{x\}\in [0,1)$ denotes the fractional part of $x$. (Note this
is slightly different from  $\lp x \rp\in (0,1]$, which appears in
the definition of $s_{i,j}$: by definition, $\lp 0 \rp=1$, while $\{0\}=0$.)
One then finds that
$$
\ft_{\sigma}(x,y)=qxy \frac1{1-e^{-qx}}\frac1{1-e^{-qy}}
\sum_{k=0}^{q-1}
e^{-q( \frac{k}{q}x +\{\frac{kp}{q}\} y )}.
$$

We may then compute $f_{i,j}(p,q)$ as the coefficient of $x^iy^j$ in the above
expression.  It is convenient to replace $x$ and $y$ with $-x$ and $-y$, which
introduces a factor of $(-1)^{i+j}$.  We obtain
$$
f_{i,j}(p,q)=(-1)^{i+j}q^{i+j-1}
\sum_{k=0}^{q-1} 
\text{coeff}\biggl(x^i; \frac{xe^{\frac{k}{q}x}}{1-e^x}\biggr)
\text{coeff}\biggl(y^i; \frac{ye^{\{\frac{kp}{q}\}y}}{1-e^y}\biggr).
$$

It is then clear that we can write the above sum in terms of values of the
Bernoulli polynomials, as follows:
 $$
f_{i,j}(p,q)=(-1)^{i+j}q^{i+j-1}
\sum_{k=0}^{q-1} 
\Ber_i\biggl(\frac{k}{q}\biggr)
\Ber_j\biggl(\{\frac{kp}{q}\}\biggr).
$$

Using the identity 
$$
\Ber_j(\lambda)=(-1)^j \Ber_j(1-\lambda),
$$
we may rewrite our expression as
 $$
f_{i,j}(p,q)=(-1)^{i}q^{i+j-1}
\sum_{k=0}^{q-1} 
\Ber_i\biggl(\frac{k}{q}\biggr)
\Ber_j\biggl(\lp -\frac{kp}{q} \rp\biggr).
$$

The sum on the right hand side is easily seen to equal the sum defining $s_{i,j}(p,q)$, 
except for a possible discrepancy in the $k=0$ term.
It follows easily from the definitions that the $k=0$ terms actually match
unless $i=j=1$, or $i=1$ and $j$ is even, or $j=1$ and $i$ is even.  In these
cases, we need the correction terms which appear in the statement of the
Theorem.


\ifx\undefined\bysame
	\newcommand{\bysame}{\leavevmode\hbox to3em{\hrulefill}\,}
\fi

\end{document}